\documentclass[11pt,conference,twocolumn]{IEEEtran}
\usepackage{amsmath,cite,latexsym,amssymb,amscd,psfrag}
\usepackage{graphicx,subfigure}

\newcommand{\covering}{\widehat{\in}}
\newcommand{\uncovering}{\widehat{\neq}}

\title{{\bf Recursive-Cube-of-Rings (RCR) Revisited: Properties and Enhancement}}
\author{{\bf Kai Xie}$^1$, \ \ {\bf Jing Li}$^1$,  \ \ {\bf Yumei Wang}$^2$, \ \ {\bf Chau Yuen}$^3$ \\
\small{1. Electrical and Computer Engineering Department, Lehigh
University, Bethlehem, 18015, USA}\\
\small{2. School of Information and Communication Engineering,} \\
\small{Beijing University of Posts and Telecommunications, Beijing, 100876 China} \\
\small{3. Singapore University of Technology and Design, Singapore, 279623}\\ 
\small{Emails: kax205@lehigh.edu, jingli@ece.lehigh.edu (contacting author), ymwang@bupt.edu.cn, yuenchau@sutd.edu.sg}
}

\begin{document}
\twocolumn[
  \begin{@twocolumnfalse}
   \maketitle

\hrule 
\vspace{0.25cm}
\small{{\bf Abstract:} We study recursive-cube-of-rings (RCR), a class of scalable graphs that can potentially provide rich inter-connection network topology for the emerging distributed and parallel computing infrastructure.  Through rigorous proof and validating examples, we have corrected previous misunderstandings on the topological properties of these graphs, including  node degree, symmetry, diameter and bisection width. 
    To fully harness the potential of structural regularity through RCR construction, new edge connecting rules are proposed. The modified graphs, referred to as {\it Class-II RCR}, are shown to possess uniform node degrees, better connectivity and better network symmetry, and hence will find better application in parallel computing.}   

\vspace{0.25cm}
\hrule
\vspace{0.35cm}
 \end{@twocolumnfalse}
]

{
  \renewcommand{\thefootnote}%
    {\fnsymbol{footnote}}
  \footnotetext[1]{Li's work is partially supported by US NSF under Grants No. CMMI-0928092, CCF-0829888, and OCI-1122027. Wang's work is partially supported by NSFC under Grant No. 61171098. Yuen's work is partially supported by the International Design Center (IDG31100102 and IDD11100101).}
}


\section{Introduction}

Cloud computing and cloud storage provide unprecedented computing, 
storage and information processing capabilities. 
Unlike last century's mainframe, today's cloud computing/storage systems
are usually comprised of hundreds of thousands of
processing elements (PE) that interconnect and process in a highly parallel, efficient and trust-worthy manner.   
These parallel computation systems may either operate on shared memory/storage or distributed memory/storage, where the latter scales better and works better for massive processing elements \cite{voTamVan2012}.

In a distributed memory/storage system, the PEs connect with and communicate to each other through an {\it interconnection network} \cite{Grammatikakis2001,Dally2004,Mesbahi2010}. Mathematically, an interconnection network is a non-directed graph with no parallel edges or self-loops, where the processing elements (e.g. personal computers) serve as the vertices and the connecting wires (e.g. optical fibers) serves as the edges. To fully harness the power provided by the distributed PEs requires the supporting interconnection network to be judicious organized with efficient communication, low hardware cost, easy applicability of algorithms, strong scalability and fault-tolerance. Although dynamic interconnection (such as switch networks) is also available, most interconnection networks use static interconnection, whose topology is critical to 
the performance and the cost of the parallel system. For example, a complete graph provides
efficient one-hop communication between any two PEs, but planar complete graphs exist for at the most 4 nodes, and a single chip with  5 or more processors must  therefore use the more complicated and expensive multi-layer design.

Researchers have developed a number of metrics to measure the goodness of a network topology, reflecting either the performance or the cost or both.
For instance, the vertex degree reflects the hardware cost, the bisection width
indicates the efficiency of the communication across the network and level of disrupt-tolerance, and the network diameter
 reveals the maximum communication delay. Additionally, network symmetry in general simplifies the resource management and provides easier means to apply 
algorithms than an asymmetric topology.

The properties of many basic network models , such as tree, ring, mesh, torus hypercube, butterfly, and de Bruijin networks,  have been well
studied (e.g.\cite{Grammatikakis2001,Dally2004,Gavriloveska2011}). Using the method of mutation or crossover, these basic models have also been modified, extended, or integrated to provide richer and better interconnection topologies. For example, modified hypercubes are proposed to
improve certain topological properties of hypercubes: Folded-hypercube
\cite{El-Amawy1991} and cross hypercube \cite{Efe1992} offer a
smaller diameter,  meta-cube \cite{Li2002} and
exchanged hypercube \cite{Chen2007} require a lower
hardware cost, and self-similar cubic \cite{Chu2012} provides a better scalability, communication and on-chip fabrication possibility. A variety of crossover constructions based on two or more basic network models are also available (e.g. \cite{Efe1996,Shi1998,Xu2010}). Among them is recursive-cube-of-rings (RCR), a family of highly-scalable crossover networks that
were originally proposed in \cite{Sun2000} and further
studied in \cite{Hu2005}.

An RCR network 
consists of well-structured cubes and rings, and is generated from
recursive expansion of the generation seed. Routing algorithms were subsequently proposed for RCR interconnection networks (e.g. \cite{choi2008}. It was shown in
\cite{Sun2000,Hu2005} that  RCR networks can possess such desirable
properties  as scalability, symmetry, uniform node degrees, low
diameter, and high bisection width. However, caution should be
exercised in choosing the parameters for RCR networks, since not all
RCR networks enjoy good properties, and many parameters will result
in asymmetric and/or unconnected networks.

The inventors defined the RCR network and analyzed some of the most
important properties of an interconnection network, including the
node degree, the bisection width and the diameter \cite{Sun2000}.
While the construction of RCR networks through recursive expansion
is rather straight-forward, RCR networks can take rather diverse
forms depending on the parameters used. The observations and
conclusions made in \cite{Sun2000} (about node degree, the bisection width and the diameter) spoke for only some RCR cases,
and did not cover all the possible scenarios. These flaws were noted
by the authors of \cite{Hu2005}, and improvements were made to the
original results. However, a careful investigation reveals that the
results provided in \cite{Hu2005} about these important properties remain incomplete.

The purpose of this paper is to rectify and improve the results in
\cite{Sun2000} and \cite{Hu2005} and to provide a more complete and
accurate characterization of RCR networks. In the first part of the
paper, we analyze the node degree (Section \ref{sec:node degree}), the bisection
width (Section \ref{sec:bisection}) and the diameter (Section \ref{sec:diameter}) of the RCR networks, and provide illustrating
examples to support our discussion.

In the second part, we further propose a class of modified RCR networks, thereafter referred to as {\it Class-II RCR} and denoted as RCR-II (Section \ref{sec:modified RCR network}). Through network analysis of connectivity, node degree, network symmetry, diameter, and bisection width, we show that the new class of RCR networks possess better structural and topological regularity  than the original RCRs defined in \cite{Sun2000}. For example, RCR-II networks have larger bisection widths and shorter network diameters than their RCR counter-parts.  An RCR-II network is guaranteed to have uniform vertex degree and, if the parameters are properly chosen, also exhibit desire symmetry property.
 The findings of this paper will help clarify and correct the misconceptions on the original RCR networks, illuminate a new and better means to exploit the structure of recursive-cube-of-rings, and provide a guideline for
choosing good parameters for RCR networks.

\section{RCR Networks}
\label{sec:RCR network}

\subsection{Introduction of RCR Networks}
\label{sec:RCR}

An RCR network consists of a host of rings connected by cube links
\cite{Sun2000}. The structure of an RCR network is completely
determined by a triple of parameters, the dimension of the cube $k$,
the size of a ring $r$, and  the level of the expansions $j$ from
the generation seed, and is thereafter denoted as RCR$(k,r,j)$.
\begin{itemize}
\item When $j=0$, we have the seed network RCR$(k,r,0)$ from which RCR$(k,r,j\ge 1)$ expands.
\item When $r=1$, the rings degenerates to a single point and the RCR network reduces to a hyper-cube. In other words, RCR networks subsume hyper-cubes as their special case.
\item When $k=0$, the cube vanishes, and the RCR$(0,r,j)$ network becomes a set of $j+1$ disconnected rings each of size $r$. Such is of little value to parallel computing.
\end{itemize}
For convenience, in the discussion that follows,  we assume
 $k\ge 1$, $r\ge 1$ and $j\ge 0$.

Let us briefly summarize the structure of an RCR$(k,r,j)$ network
\cite{Sun2000} and introduce the notations that will be used in the
discussion.

\begin{figure*}[htbp]
\centering
\includegraphics[width=4in]{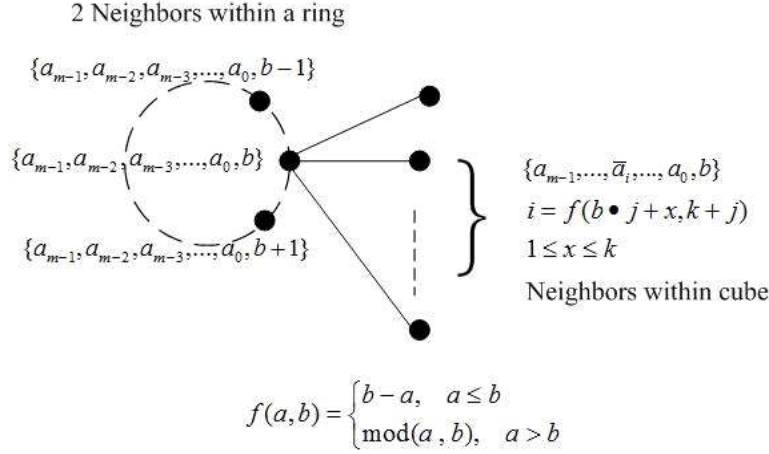}
\caption{Construction of an RCR$(k,r,j)$ network} 
\label{fig:neighbors}
\end{figure*}

An RCR$(k,r,j)$ network has altogether $2^{k+j}r$ nodes in the network.
As shown in Figure \ref{fig:neighbors},
 each node in RCR$(k,r,j)$ is represented by its coordinate,
$<a_{k+j-1},a_{k+j-2},...,a_{0};b >$, which consists of a cube
coordinate $<a_{k+j-1},a_{k+j-2},...,a_{0}>$ and a ring coordinate $b$.
An RCR$(k,r,j)$ network is expanded from  RCR$(k,r,j-1)$ by replicating RCR$(k,r,j-1)$ twice, preserving all the ring edges, and breaking and reconnecting all the cube edges.

The ring coordinate $b$, where $0\leq b \leq r-1$, specifies the position of the node within a ring of dimension $r$. When $r=1$ and $2$, the ring reduces to a single node and a single line, respectively. When $r>2$, each node has two
distinct ring neighbors that having the same cube coordinate but adjacent ring coordinates,
\begin{align}
& <a_{k+j-1},a_{k+j-2},...,a_{0}; b\stackrel{.}{-}1>, \\ 
\rm{and} \ &
<a_{k+j-1},a_{k+j-2},...,a_{0}; b\stackrel{.}{+}1>,
\end{align}
where $\stackrel{.}{+}$ and $\stackrel{.}{-}$ stand for the modulo $r$
arithmetic:
\begin{align}
\alpha \stackrel{.}{\pm}\beta= \mod(\alpha \pm \beta, r).
\end{align}
In general, a node has $\min(r-1,2)$ ring neighbors.

The cube coordinate $<a_{k+j-1},a_{k+j-2},...,a_{0}>$ consists of
$k+j$ binary cube bits $a_i \in\{0,1\}$. Let  $\bar{\alpha}$ denote the binary complementary of $\alpha$, such that $\bar{0}=1$ and $\bar{1}=0$.
A cube link can only exist between node $<a_{k+j-1},...,a_{i+1},a_i,a_{i-1},...,a_0;b>$ and node $<a_{k+j-1},...,a_{i+1},\bar{a_i},a_{i-1},...,a_0;b>$, and it exists only when the index $i$ and the ring coordinate $b$ satisfy the constraint
\begin{equation}
i= f(b\times j+x,k+j),\label{equ:RCR_f}
\end{equation}
where $0 \leq b \leq r-1$, $1 \leq x \leq k$, and function $f$ is defined as \cite{Sun2000}
\begin{equation}
f(a,b) = \left\{ \begin{array}{cc}
        b-a, & a \le b,\\
        \mod(a, b), & a > b
    \end{array}\right. 
\label{equ:definition_f}
\end{equation}
where $a$ and $b$ are integers. 
Hence, the number of cube neighbors which a node has is  the number
of distinct values  $f(b\times j+x,k+j)$ may take. It is easy to see
that $f(b\times j+x,k+j)$ may take at the most $k$ distinct values.

Having introduced the RCR network as defined in \cite{Sun2000}, below
we discuss the properties of RCR networks.

\section{Node Degree}
\label{sec:node degree}

The degree of a node, denoted as $D_n$, is defined as the total
number of distinct ring neighbors and cube neighbors this node has.
For an RCR network, $D_n\le k+2$. It was claimed in \cite{Sun2000}
and \cite{Hu2005} that an RCR$(k,r,j)$ network always had a uniform
node degree $D_n$ and  a symmetric structure regardless of the
parameters used. It was shown that, when the dimension of the rings
$r$ was less than or equal to 2,  the node degree was $D_n=k+r-1$
unanimously, and when $r$ was larger than $2$, the node degree
became $D_n=k+2$ unanimously \cite{Sun2000,Hu2005}.

However, this conclusion is based on the assumption that
 each node in an RCR$(k,r,j)$ network has $\min(r-1,2)$ distinct ring neighbors, and $k$ distinct cube neighbors due to the $k$ possible values of
$x$ in (\ref{equ:RCR_f}), where $1 \leq x \leq k$. This assumption is valid for all the RCR examples presented in \cite{Sun2000}, but does not hold in general.
Depending on the choice of the  cube dimension
 $k$, ring dimension $r$ and the expansion level $j$, different values of $x$ may generate the same value of $f(b\times j +x, k+j)$ (here the auxiliary variable $b$ is an integer, $0 \leq b \leq
r-1$). It is therefore possible for a node to have fewer than $k$ distinct cube neighbors. 

\begin{figure*}
\centerline{
\includegraphics[width=4.0in]{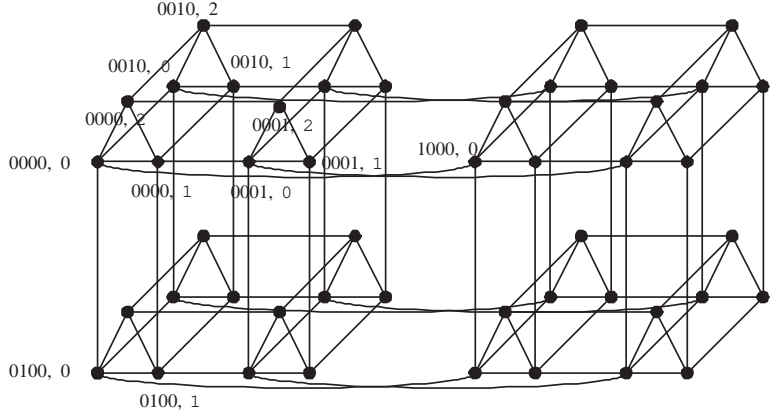}
}
\caption{Construction of RCR(3,3,1) network: nonuniform-node-degree
RCR network}
 \label{fig:nodedegree}
\end{figure*}

\vspace{0.2cm}
{\bf \emph{ Example 1: [Non-uniform node degree of RCR]}} Consider an
RCR$(3,3,1)$ network as shown in fig. \ref{fig:nodedegree}.
Following the definition in Subsection \ref{sec:RCR} and in
\cite{Sun2000}, $x$ may take  three possible values:
$x\in\{1,2,3\}$, and the auxiliary variable $b$ may also take  three
possible values: $b\in\{0,1,2\}$. A node in this RCR network may
experience one of the following three scenarios:
\begin{itemize}
  \item When $b=0$, the possible values for $f(bj +x, k+j)$
  are $3,2,1$, which correspond to  $x=1,2,3$, respectively. The nodes in this case will have 3 cube neighbors, which lead to a node degree of $D_n=3+2=5$.
  \item When $b=1$, the possible values for $f(b j +x, k+j)$
  are $2,1,0$, so the nodes here also have 3 cube neighbors and a node degree of 5.
  \item When $b=2$, there are 2 possible values for $f(b j +x, k+j)$: with $x=1$ and $x=3$ we have $f(2\cdot 1 + 1, 4) =1=f(2\cdot 1+3, 4)$, and with $x=2$ we have $f(2 \cdot 1 + 2, 4)=f(4,4)=0$. The nodes here thus have only 2 cube neighbors and a degree of $4$.
\end{itemize}
Since the node degrees are not uniform, the RCR(3,3,1) network cannot be symmetric.

\vspace{0.2cm}
{\bf \emph{Theorem 1: [Node degree of RCR]}} The node degree $D_n$ of an RCR$(k,r,j)$ network satisfies
\begin{itemize}
  \item When $r \leq 2$, $D_n = k+r-1$ for all the nodes;
  \item When $r > 2$ and at least one of $k \leq j+1$ or $j = 0$ is satisfied,  $D_n = k+2$ for all the nodes;
  \item When $r >2$, $j \geq 1$ and $k >j+1$, $D_n$ is not a constant but takes multiple values:    $\lceil k/2 \rceil +2 \leq D_n \leq k+2$.

\end{itemize}

{\it Proof:} Case I: $r \leq 2$. Each node has $r-1$ ring neighbors.
Since $b j + x  \leq  (r-1) j + k  \leq  j+k$,
following the definition of $f(\cdot)$ in (\ref{equ:definition_f}), we have $f(b j + x, k+j) = (k+j)-(bj+x)=k+(1-b)  j-x$. Consider a
node with given parameters $b$,$j$,$k$. The function $f(\cdot)$ is a linear
function of $x$ and generates $k$ distinct output values for $k$ distinct input values $x$, indicating that a node always has $k$ cube neighbors. The node degree is therefore $D_n=k+r-1$.

Case II: $r > 2$ and $k \leq j+1$.  Each node here has 2 ring neighbors.
To evaluate the number of cube neighbors, consider separating
 the nodes in two cases: $b \leq 1$ and $b> 1$.
(i) When $b \leq 1$, we have $b j + x <j+k$ and hence $f(b
j + x, k+j) = k+ (1-b)  j-x $, which assumes $k$ distinct
values for $x=1,2,\cdots,k$. (ii) When $b>1$, we have $ b j + x  \geq   2  j + 1 \geq j+k$, and therefore $f(b j + x, k+j) = mod(b  j +x, k+j)$, which again generates $k$ distinct values with input $x\in[1,k]$. In either case, a node has $k$ cube neighbors and 2 ring neighbors, making the node degree $D_n=k+2$.

Case III: $r >2$ and $j = 0$. Each node has 2 ring neighbors. Since
$j=0$, the possible values of $f(b j + x,k)=f(x,k)=k-x$ are
$\{k-1,k-2,...,0\}$ for $1 \leq x \leq k$. There exist $k$ distinct
values for $f(b j + x,k)$. Therefore, the node degree should be
$D_n=k+2$.

Case IV: $r > 2$, $j \geq 1$ and $k > j+1$. Each node here has 2
ring neighbors. Separate all the nodes in three cases: $b \leq 1$,
$b = 2$, or $b>2$. (i) When $b \leq 1$, $bj+x\leq j+x\leq j+k$ and
$f(b  j + x, k+j) = k+ (1-b)  j-x$, yielding $k$ distinct values. So
each node has $k$ cube neighbors and a degree of $D_n=k+2$.
 (ii) When $b = 2$, since $k-j > 1$, we have
$  b j + k-j = j +k$. There always exists a positive integer $t\leq
k-j-1$ and $t \leq j$. Then we have $1 \leq k-j-t < k-j+t \leq k$
and $b j + (k-j-t) < k+j < b j + (k-j+t)$. Let $x_1=k-j-t$ and
$x_2=k-j+t$. Thus,
\begin{align}
  f(bj+x_1)&= f(2j + k-j-t, k+j),\nonumber\\
     & =   (k+j)-(j+k-t)= t,\label{equ:theorem1_3}\\
  f(bj + x_2, k+j) &= f(2j+ k-j+t,k+j),\nonumber\\
  &=mod(k+j+t ,k+j)= t.
\end{align}
Since $ x_1 \neq x_2$, the nodes here have fewer than $k$ cube neighbors, and hence their node degree is strictly smaller
 than $k+2$. Comparing (i) and (ii), we
know that the node degree can not be uniform in Case III. (iii)
Following the same procedure, we can shown that when $b>2$,  the
node degree may be either equal to or smaller than $k+2$. In
conclusion, when $r > 2$, $j \geq 1$ and $k > j+1$, the node degree
is not fixed.

From the above discussion, we have known that, for given $b>1$, the
node degree may be less than $k+2$. This comes from the fact that
some possible values of $bj+x$ are less than $k+j$ and the others
are larger than $k+j$. Some $bj+x$ less than $k+j$ will give the
same $f(bj+x, k+j)$ with certain $bj+x$ larger than $k+j$.
Obviously, the overlapping part is at most $\lfloor k/2 \rfloor$.
Therefore, the node degree is always not less than $\lceil k/2
\rceil +2$.

\vspace{0.2cm}
{\bf \emph{Remark:}} It should be noted that a uniform node degree is but a
necessary condition for a network to be symmetric. Network symmetry
is a stronger condition than merely having a uniform node degree and
some apparent regularity in structure. In the case of RCR networks,
despite their well-defined structure, a uniform node degree does not
necessarily lead to network  symmetry.

\begin{figure*}[htbp]
\centering
\includegraphics[width=4.0in]{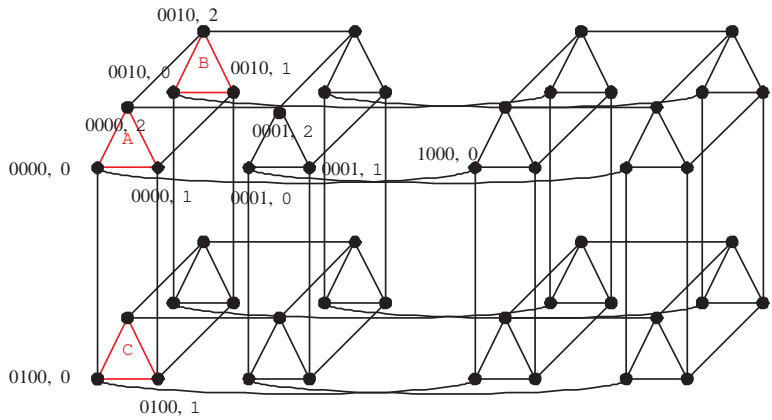}
\caption{RCR(2,3,2): A uniform-node-degree but asymmetric network.}
\label{fig:symmetric}
\end{figure*}

\vspace{0.2cm}
{\bf \emph{Example 2: [Uniform-node-degree but asymmetric RCR]}} The
RCR(2,3,2) network shown in Figure \ref{fig:symmetric} has a uniform
node degree $D_n=4$, but is asymmetric. To see this, consider
setting an arbitrary node in RCR(2,3,2)  as $<0000;0>$ and
relabeling all the nodes. From the definition of symmetry, the
newly-labeled network will preserve the same connection as the
original one. Suppose we relabel node $s=<0000;1>$ as $<0000;0>$.
Observe that node $s$ has two neighboring rings: ring $B$ and ring
$C$, and both rings have two cube edges connected with ring $A$
which node $s$ belongs to. It is impossible to find a relabeling
scheme for the same network structure that will satisfy the rules
(definition) of RCR networks.

\section{Bisection Width}
\label{sec:bisection}

The bisection width, defined as the minimum number of edges that
must be removed in order to bisect a network, is another important
property for interconnected networks. The conclusion made in
\cite{Sun2000} on the bisection width of RCR networks was not
entirely correct. Some missing cases were picked up and amended in
\cite{Hu2005}, but others remain overlooked.

In \cite{Sun2000}, the bisection width, $B_{RCR}(k,r,j)$, is computed as
\begin{equation}
\label{eqn:width}
B_{RCR}(k,r,j)=
Num(k,r,j)\times N/(2\times r),
\end{equation}
for all ring dimensions $r$, where $N$ is the total number of nodes
in the networks, and $Num(k,r,j)$ is defined as the number of $b$
values satisfying $f(b j + x, k+j) = k+j-1$.

The authors of \cite{Hu2005} recognized that the case of $r=1$ is an
exception. They showed that an RCR$(\cdot,1,\cdot)$ network
comprises two  unconnected subnetworks of equal sizes and therefore
has 0 bisection width. They thus amended the results in
\cite{Sun2000} by setting condition $r\ge 2$ on (\ref{eqn:width}),
and adding the case of $B_{RCR}(k,1,j)=0$. Below we show that more
exceptions exist such that an RCR network with $r\ge 2$ may still be
unconnected and has 0 bisection width.

\begin{figure*}[htbp]
\centering
\includegraphics[width=4.5in]{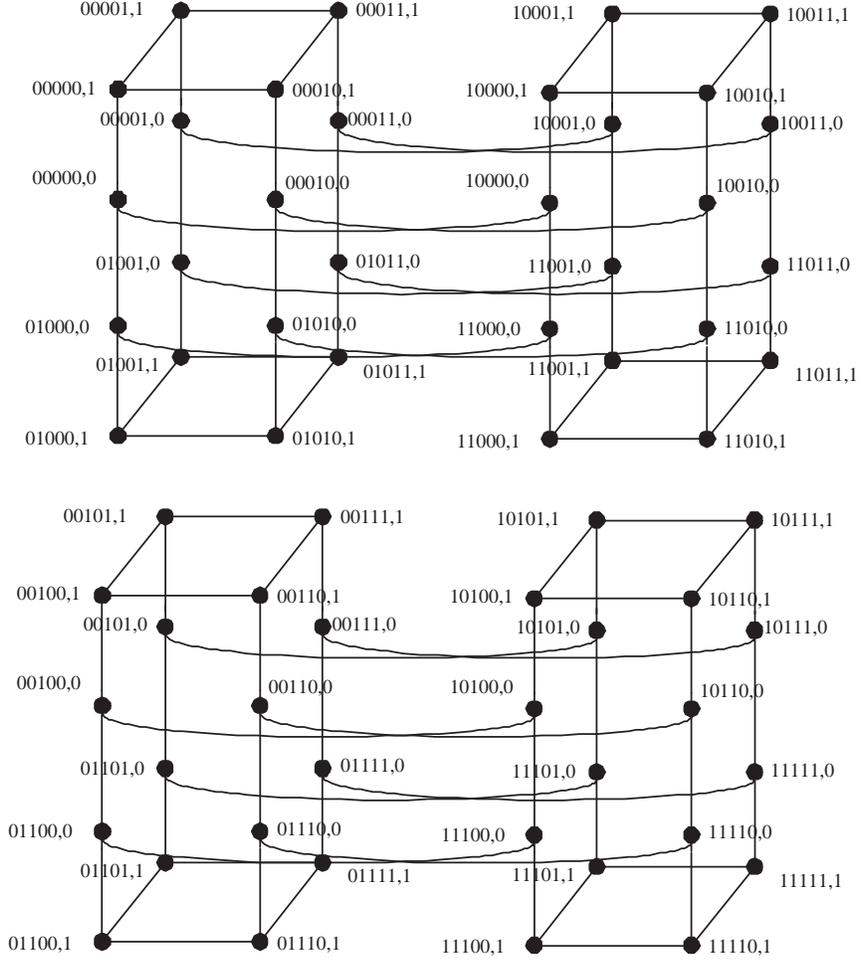}
\caption{Construction of RCR(2,3,2) network: unconnected RCR network
with $r=2$} 
\label{fig:unconnected}
\vspace{-0.3cm}
\end{figure*}

{\bf \emph{Example 3: [An unconnected RCR network with $r=2$]}} Consider the
RCR(2,2,3) network in Figure \ref{fig:unconnected},  whose nodes
have coordinates $<a_4, a_3, a_2,a_1, a_0; b>$. When $b=0$, the
possible values of $f(b j + x, k+j)$ are $4$ and $3$; when $b=1$,
the possible values of $f(b j + x, k+j)$ are $0$ and $1$. In other
words, a node can only have a cube neighbor whose coordinate differs
from that itself in one of the four bit positions $a_4$, $a_3$,
$a_1$ and $a_0$. Thus the two sets of nodes, $\{<a_4, a_3, 0,a_1,
a_0; b>\}$ and  $\{<a_4, a_3, 1,a_1, a_0;b>\}$, each consisting of
$2^4\cdot 2=32 $ nodes, do not have any inter-connecting edge
between them. The network is thus unconnected and has a bisection
width of 0.

The reason that \cite{Sun2000,Hu2005} failed to spot such cases as Example 3 is that, in computing the bisection width, they always bisected the network into two sub-networks with the $(m-1)$th bit being the complementary of each other. However, this cut is not always the minimum cut. To help evaluate  the bisection width,  let us introduce a new parameter $Num(k,r,j,t)$.

\vspace{0.2cm}
{\bf \emph{Definition 1: [$Num(k,r,j,t)$]}} Consider an RCR$(k,r,j)$ network. For a given integer $b$
if there exists an integer $x\in[1,k]$ such that $f(b j + x, k+j) = t$,
we say $b$ satisfies $f(b j + x, k+j) = t$. $Num(k,r,j,t)$ is defined as
the number of integer values $b\in [0,k-1]$ that satisfies $f(b j + x, k+j) = t$.

\vspace{0.2cm}
{\bf \emph{Theorem 2: [Bisection Width of RCR]}} The bisection width of an RCR$(k,r,j)$ network is upper-bounded by:
\begin{equation}
B_{RCR}(k,r,j)\le \min_{t\in\{0,\cdots,k+j-1\}} (Num(k,r,j,t))\times N/(2\times r),
\end{equation}
where $r$ is the dimension of rings and $N$ is the total number of nodes.

{\it Proof}: Consider bisecting the network nodes into two groups,
$\{<a_{k+j-1},...,a_{t+1},a_t=0,a_{t-1}...,a_0;b>:  a_i\in\{0,1\} \forall i \rm{\ except\ } i\ne t, b\in\{0,1,\cdots, r-1\}$ and
$\{<a_{k+j-1},...,a_{t+1},\bar{a}_t=1,a_{t-1},...,a_0;b>: a_i\in\{0,1\} \forall i \rm{\ except\ } i\ne t, b\in\{0,1,\cdots, r-1\}$, where $0\le t\le k+j-1$. A cube edge exists between node $<a_{k+j-1},...,a_{t+1},0,a_{t-1}...,a_0;b>$ and node $<a_{k+j-1},...,a_{t+1},\bar{a}_t=1,a_{t-1},...,a_0;b>$ if and only if $b$ satisfies $f(b\times j + x, k+j) = t$. Following the definition, there are  $Num(k,r,j,t)$ different values of $b$ satisfying $f(b\times j + x, k+j) = t$. Given $t$ and $b$, there exist $2^{k+j-1} = N/(2r)$ possible values for $a_{k+j-1}, \cdots, a_{t+1}, a_{t-1},\cdots, a_0$. Therefore, the are altogether $Num(k,r,j,t)\times N/(2 r)$ edges between the two groups. Hence, the bisection width is
$\min_t (Num(k,r,j,t))\times N/(2\times r)$, where  $0\leq t \leq k+j-1$.

\vspace{0.2cm}
{\bf \emph{Remark:}} Theorem 2 considers the case where a bisection cut
consists of cube edges only. It is possible for a set of ring edges
to also form a bisection cut and to have a smaller size than those
formed from cube edges. Hence, what is provided in Theorem 2
represents an upper bound rather than the exact bisection width, as
shown in example 4. However, this upper bound is tight, as shown in
example 5.

\begin{figure}[htbp]
\centering
\includegraphics[width=3.0in]{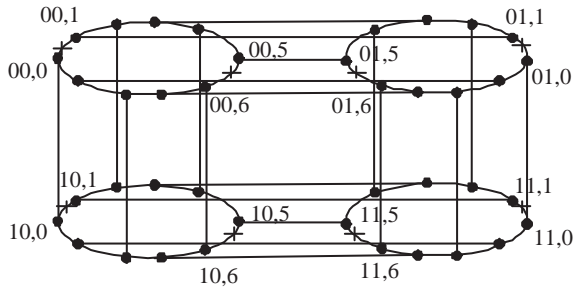}
\caption{Construction of RCR(1,10,1) network: Upper bound is not
exact bisection width} \label{fig:example4}
\end{figure}

\vspace{0.2cm}
{\bf \emph{Example 4: [Bisection cut may be formed by ring edges]}}
  Consider
an RCR$(1,10,1)$ network, which comprises 4 rings of dimension 10
each, 10 cube edges connecting node pairs $<c0;b>$ and $<c1;b>$ for
$b=1,3,...,9$ and $c=\{0,1\}$; and another 10 cube edges
connecting node pairs $<0c;b>$ and $<1c;b>$ for $b=0,2,...,8$ and
$c=\{0,1\}$.
To bisect the network through cube edges, the minimum cut consists
of 10 cube edges. However, the minimum bisect width is 8, resulted
from 8 ring edges that connect, say, $<00;0>$ and $<00;1>$, $<00;5>$
and $<00;6>$, $<01;5>$ and $<01;6>$, $<01;0>$ and $<01;1>$, and
$<10;0>$ and $<10;1>$, $<10;5>$ and $<10;6>$, $<11;5>$ and $<11;6>$,
$<11;0>$ and $<11;1>$, as shown in Fig. \ref{fig:example4}.

\begin{figure}[htbp]
\centering
\includegraphics[width=2in]{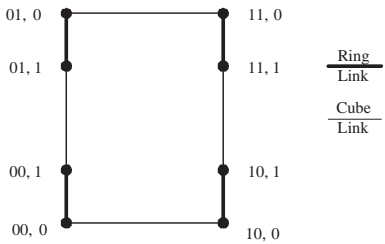}
\caption{Construction of RCR(1,2,1) network: upper bound is tight}
\label{fig:bisection_bound}
\end{figure}

\vspace{0.2cm}
{\bf \emph{Example 5: [Upper bound of bisection width is tight]}}
Consider an RCR$(1,2,1)$ network as shown in
Fig.\ref{fig:bisection_bound}. According to Theorem 2,
$\min_{t\in\{0,\cdots,k+j-1\}} (Num(k,r,j,t))=1$ and $B_{RCR}(k,r,j)
\leq 2$. From Fig.\ref{fig:bisection_bound}, it is easy to see that
the bisection width is exactly 2, which achieves the bound in Theorem 2 with equality.

By definition, an interconnected network is not supposed to be unconnected. The fact that an RCR network may be unconnected (see Example 3) suggests that one needs to exercise with caution in choosing the parameters. Below we present the necessary and
sufficient condition that guarantees the connectivity of an RCR network.
 For convenience, we introduce  a new notation $\covering$.

\vspace{0.2cm}
{\bf \emph{Definition: 2}} If $A$ can take all the possible integer values
between $\alpha$ and $\beta$ (inclusive), we say that $A$ {\it
covers} the range $[\alpha, \beta]$, and denote it as
$A\covering[\alpha, \beta]$. Otherwise, we say that range $[\alpha,
\beta]$ is not covered by $A$ and denote it as $A\uncovering[\alpha,
\beta]$

\vspace{0.2cm}
{\bf \emph{Theorem 3: [Sufficient and necessary condition for RCR to be
connected]}} An RCR$(k,r,j)$ network is connected if and only if
 $f(b j + x, k+j)\covering [0, k+j-1]$ for $0\le b\le r-1$ and $0\le x\le k$.

{\it Proof:} (Sufficient condition: when $f(bj+x,k+j)$ covers all
the integer values between 0 and $k+j-1$, then the RCR network is
connected.) Notice that all the nodes having the same cube
coordinates (but different ring coordinates) form a ring, and hence
we can use the common cube coordinate to identify a ring. To show
the connectivity of an RCR$(k,r,j)$ network, it is sufficient to
show that any two ``adjacent'' rings are connected, where by
adjacent, we mean that the cube coordinates of the two rings differ
in only one bit position. Consider two adj rings with respective
cube coordinates $<a_{k+j-1} \cdots a_{t_0+1}0 a_{t_0-1}\cdots a_0>$
and  $<a_{k+j-1} \cdots a_{t_0+1}1 a_{t_0-1}\cdots a_0>$. Since
$(bj+x,k+j)\covering [0,k+j-1]$, $f(b_0 j+x_0,k+j)=t_0$ for some
valid values of $b_0$ and $x_0$. According to the definition of RCR
networks, a cube edge exists that connects node  $<a_{k+j-1} \cdots
a_{t_0+1}0 a_{t_0-1}\cdots a_0;b_0>$ and node $<a_{k+j-1} \cdots
a_{t_0+1}1 a_{t_0-1}\cdots a_0;b_0>$.    Hence any source node in
the first ring can travel through ring edge(s) to reach node
$<a_{k+j-1} \cdots a_{t_0+1}0 a_{t_0-1}\cdots a_0;b_0>$, and then
through the cube edge to get to $<a_{k+j-1} \cdots a_{t_0+1}1
a_{t_0-1}\cdots a_0;b_0>$ in the second ring, and again through ring
edge(s) to get any destination node in the second ring.

(Necessary condition: when the RCR network is connected, then
$f(bj+x,k+j)$ covers all the integer values between 0 and $k+j-1$.)
Proof by contradiction. Suppose that the network is connected but
there exists $t_0\in[0,k+j-1]$ such that $f(bj+x,k+j)\ne t_0$ for
all the valid values of $b$ and $x$. According to the definition of
RCR  networks, there does not exist a cube edge connecting any pair
of nodes $<a_{k+j-1} \cdots a_{t_0+1}0a_{t_0-1}\cdots a_0; b>$ and
$<a_{k+j-1} \cdots a_{t_0+1}0a_{t_0-1}\cdots a_0; b>$. Hence, an
arbitrary node whose coordinate has the former form  and an
arbitrary node whose coordinate has the latter form are not
reachable to each other. For example, node $<0\cdots 0;0>$ cannot
reach node $<1\cdots 1; 0>$. This contradicts with the connectivity
assumption.

\vspace{0.2cm}
{\bf \emph{Lemma 4: [Necessary condition for  RCR to be connected]}} An
RCR$(k,r,j)$ network is unconnected, if $(r-1) k < j$.

{\it Proof}: A node in an RCR$(k,r,j)$ network is denoted by the
combination of a cube coordinate and a ring coordinate, where the
former is a length-$(k+j)$ binary vector, and the latter takes $r$
possible values. According to the cube-edge connecting rule in
(\ref{equ:definition_f}), there are $k$ possible values for $x$ and
$r$ possible values for $b$, and hence at the most $kr$ possible
values for  $f(b j + x, k+j)$.  When $(r-1)k <j$, or, $  rk < k+j$,
there must be at least one bit index $t$ in the length-$(k+j)$ cube
coordinate that does not equal any value of $f(b j + x, k+j)$.
Hence,  $f(b j + x, k+j)\uncovering [0, k+j-1]$. According to
theorem 3, the network is therefore unconnected.

\vspace{0.2cm}
{\bf \emph{Theorem 5: [Parameters for RCR to be connected]}}
An RCR$(k,r,j)$
network is connected, if and only if
\begin{equation*}
    \label{equ:definition of f}
     \left\{ \begin{array}{ll}
        (r-1)k \geq j, & r \leq 2\\
        (r-1)k \geq j+1, & r > 2
    \end{array}\right.
\end{equation*}

{\it Proof}: From Theorem 3, it is sufficient to show that these
parameters, and only these parameters, ensure that function
$f(bj+x,k+j)\covering [0, k+j-1]$.

Case I: $r=1$. When $r=1$ and $(r-1)k \geq j$, then $j=0$, and the
only valid value for $b\in[0,r-1]$ is 0. Thus, $f(b j + x,
k+j)=f(x,k)=k-x \covering[0,k-1]$ for $x\in[1,k]$. It follows from Lemma 4 that  $(r-1)k
\geq j$ is the necessary condition for  $r=1$.

Case II: $r=2$. When $r=2$ and $(r-1) k \geq j$, then $j \leq k$,
and $b$ may take values of either 0 or 1.
\begin{align}
{\rm when\ } b=0,\ \ &  f(b j + x, k+j)=f(x,k+j),\nonumber \\
                    &=k+j-x\covering[j, k+j-1],\\
{\rm when\ } b=1,\ \ &  f(b j + x, k+j)=f(x+j,k+j),\nonumber\\
                    &=k-x\covering[0,k-1].
\end{align}
Since $j\le k$,  $f(b j + x, k+j)\in[0,k+j-1]$. It follows from Lemma 4 that $(r-1)k \geq
j$ is the necessary condition for $r=2$.

Case III: $r> 2$. We first show that  $(r-1)k \geq j+1$ is a sufficient condition for $f(bj+x,k+j)\covering[0,k+j-1]$ by differentiating two subcases.

(i) Suppose $k > j$. Then $b$ may take values of $0,1,\cdots, r-1$.
\begin{align}
{\rm when\ } b\!=\!0,\ & f(b j \!+\! x, k\!+\!j)\!=\!(k\!+\!j)\!-\!x\covering[j, j\!+\!k\!-\!1], \\
{\rm when\ } b\!=\!1,\ & f(bj \!+\! x, k\!+\!j)\!=\!(k\!+\!j)\!-\!(j\!+\!x),\nonumber\\
&=\!k\!-\!x\covering [0, k\!-\!1].
\end{align}
 Since $k>j$, we can see $f(bj+x,k+j)\covering [0, k+j-1]$.

(ii) Suppose $k \leq j$. From the condition $(r-1)k\geq j+1$, we get 
$(j+1)/k\le r-1$. Since $j,k,r$ are integers, we have $\lfloor
j/k\rfloor+1\le r-1$. Now $b$ can take values from $0$ to $r-1$. We
show that as $b=0,1,\cdots, \lfloor j/k\rfloor +1$, the
$f(bj+x,k=j)\covering[0,k+j-1]$. Since $k\leq j+1$ and $x\in[1,k]$,
we have
\begin{align}
& {\rm when\ } b\!=\!1, \\
&\ \  f(bj \!+\! x, k\!+\!j)\!=\!(k\!+\!j)\!-\!(j\!+\!x)=k\!-\!x\ \covering [0, k\!-\!1], \nonumber\\
&{\rm when\ } b\!=\!0, \\
& \ \  f(b j\!+\!x , k\!+\!j)\!=\!(k\!+\!j)\!-\!x\ \covering [j, j\!+\!k\!-\!1],\nonumber\\
&{\rm when\ } b\!=\!2, \\
& \ \   f(bj\! +\! x, k\!+\!j)\!=\!mod(2j\!+\!x,k\!+\!j)\covering [j\!-\!k\!+\!1, j], \nonumber\\
&{\rm when\ } b\!=\!3, \\
& \ \   f(bj \!+\! x, k\!+\!j)\!=\!mod(3j\!+\!x,k\!+\!j)
\covering [j\!-\!2k\!+\!1,j\!-\!k],\nonumber\\
 & \ \ \ \ \ \ \ \ \ \ \ \ \ \cdots  \cdots \nonumber\\
&{\rm when\ } b\!=\!\left\lfloor \frac{j}{k}\right\rfloor \!+\!1, \\
& \ \    f(bj + x, k+j) = mod(  \left\lfloor \frac{j}{k}\right\rfloor j \!+\!j\!+\!x,k\!+\!j),\nonumber\\
& \ \ \ \ \ \ \ \ \ \ \ \ \ \ \ \ \ \ \ \ \ \ \ \   \covering \left[j\!-\!\left\lfloor \frac{j}{k}\right\rfloor k\! +\!1  ,j\!-\!\left\lfloor \frac{j}{k}
\right\rfloor k  \!+\!k \right].\nonumber
\end{align}
Since $j-\lfloor \frac{j}{k}\rfloor k +1 =mod(j, k)+1\le (k-1)+1=k$,
all the integer segments in the above connect and cover the entire range of $[0, j+k-1]$.

We now show that $(r-1) k \geq j+1$ is also a necessary condition for
$f(bj+x,k+j)\covering [0, k+j-1]$, by showing that the function $f$ fails to cover $[0,k+j-1]$ otherwise. Again, we evaluate two separate cases:

(i) If $(r-1) k < j$, according to Lemma 4, the RCR$(k,r,j)$ network
is unconnected.

(ii) If $(r-1) k = j$,
\begin{align}
&{\rm when\ } b=0, \ f(bj\!+\!x, k\!+\!j)=f(x,k\!+\!j),\nonumber\\
& \ \ \ \ \ \ \ \ \ \ \ \ \ \ \ \ \ \ =k\!+\!j\!-\!x\ge j =(r\!-\!1)k>k,
\label{equ:connectivity_step1}\\
&{\rm when\ } b=1,\  f(b j \!+\! x, k\!+\!j)=f(j\!+\!x,k\!+\!j),\nonumber\\
& \ \ \ \ \ \ \ \ \ \ \ \ \ \ \ \ \ \ = (k\!+\!j)\!-\!(j\!+\!x)\!=\!k\!-\!x < k, \\
&{\rm when\ } 2\! \leq \!b \!\leq\! r\!-\!\!1,  \ f(b  j \!+\! x, k\!+\!j)\!=\!f(bk(r\!-\!\!1)\!+\!x,rk),\nonumber\\
& \ \ \ \ \ \ \ \ \ \ \ \ \ \ \ \ \ \ =mod(bk(r\!-\!1)\!+\!x,rk)\ne k,
\label{equ:connectivity_step3}
\end{align}
where the first equality in (\ref{equ:connectivity_step3}) comes from the assumption $(r-1)k=j$, and the second equality comes from the definition of function $f$. To see that the last inequality in  (\ref{equ:connectivity_step3}) holds, we use proof by contradiction: Suppose there exists $b_0\in[2, r-1]$ and $x_0\in[1,k]$ such that $mod(bk(r-1)+x,rk)=k$. That is, we can find an integer $B$ satisfying  $bk(r-1)+x=rkB+k$. Since $x$ must be an integer multiple of $k$ in order for the equality to hold, we have $x=k$. The equality now transfers to $b(r-1)+1=rB+1$, or, $br-b=rB$. Clearly, $b$ must be an integer multiple of $r$ in order for the equality to hold, but $b\in[2, r-1]$, resulting in a conflict.
Hence, it follows from (\ref{equ:connectivity_step1})-(\ref{equ:connectivity_step3}) that when $(r-1) k = j$, $f(bj+x,k+j)$ does not produce an output $k$ and hence does not cover $[0,k+j-1]$.

\vspace{0.2cm}
{\bf \emph{Corollary 6}}: If the node degree of an RCR$(k,r,j)$ network is
non-uniform,  then this RCR network is connected.

{\it Proof}: From Theorem 1, the node degree of an RCR$(k,r,j)$
network is non-uniform if and only if $r>2$ and $k>j+1$. This leads
to $r>2$ and $(r-1)k\ge j+1$, and according to Theorem 5, the
network is connected.

\vspace{0.2cm}
{\bf \emph{Remark:}} Although RCR networks have well-defined and systematic construction, their structure regularity has not been most desirable. From the analysis we performed thus far,
(i) An RCR network does not always have uniform node degree.
(ii) Even when an RCR network has a uniform degree, it is not necessarily symmetric.
(iii) An RCR network having uniform node degree may be unconnected.  An RCR network having non-uniform node degree, on the other hand, is always connected.

\section{Network Diameter}
\label{sec:diameter}

The diameter of a network measures the minimum number of hops it
takes to reach from any node to any other node in the network.
\cite{Sun2000} stated that the diameter was upper bounded by
$k+j-1+\lceil (r-1)/2 \rceil$, which
failed to differentiate between connected
and unconnected cases. \cite{Hu2005} improved the accuracy of the
results by recognizing that an RCR network is unconnected and hence
has an infinite diameter when the ring dimension $r$ is 1 (assuming
the expansion level $j>0$).  When $r>1$, \cite{Hu2005} stated that
the diameter was upper bounded by $k+j+1+\lfloor r/2 \rfloor$, which,
in fact equals $k+j-1+\lceil (r-1)/2 \rceil +2$ for any integer
value of $r$.

%

However, the results in \cite{Hu2005} have not been accurate either.
As we have shown in Theorem 5, when  $r\!>\!1$, it is also possible for
an RCR network to become unconnected and to have an infinite
diameter. Further, even in the connected case, the upper-bound
provided in \cite{Hu2005} is on the optimistic side. The proof in
\cite{Hu2005} followed the argument that one could always take a
one-hop walk from one node to its ``cube neighbor'' whose cube
coordinate differed in one bit position from that of itself and
whose ring coordinate was the same as that of itself.
Thus, after at the most $k+j+1$ hops along the cube edges,
\cite{Hu2005} decided that the source node must have reached an
intermediate node that had the same cube coordinate as the
destination node. It then took at the most  $\lfloor r/2 \rfloor$
hops, or, half the ring dimension, along the ring to reach the
destination. This argument is flawed because a pair of nodes having
the same ring coordinates and differing in one bit in cube
coordinates are not necessarily connected directly. From the
definition of RCR networks, a connecting  edge exists between two
such nodes only when their common ring coordinate $b$ meets the
constraint in \ref{equ:definition_f}. Here is a counter-example to
the conclusion drawn in \cite{Hu2005}.

\vspace{0.2cm}
{\bf \emph{Example 6: [Diameter of RCR]}} Consider an RCR$(2,5,7)$
network whose nodes are specified by $<a_8,a_7,\cdots, a_0; b>$ and
whose possible values of $f(b j +x, k+j)$  are listed in  Table
\ref{tab:RCR(2,5,7)}.
\begin{table}[htf]
  \centering
  \begin{tabular}{|c|c|c|c|c|c|}
    \hline
      & $b=0$ & $b=1$ & $b=2$ & $b=3$ & $b=4$ \\
      \hline
    $f(.)=$ & $\{8,7\}$ &  $\{1,0\}$ &  $\{6,7\}$ &  $\{4,5\}$ &  $\{2,3\}$ \\
    \hline
  \end{tabular}
  \caption{The possible values of $f(b j +x, k+j)$ for RCR(2,5,7)}
  \label{tab:RCR(2,5,7)}
\end{table}
Suppose we want to find the distance between a source node $A=<00\cdots 00;0>$ to a destination node
$B=<11\cdots 11;2>$. From Table \ref{tab:RCR(2,5,7)}, we see
that each value of ring coordinate $b$ allows the ``flip'' of only two bits in the cube coordinate. Hence all the possible values of $b$ need to be traversed in order for node $A$ to get to node $B$. The
shortest path is found as follows:

With the ring coordinate $b=0$, we flip bits $a_8$ and $a_7$, i.e. take two hops across the cube edges:
\begin{align}
<000000000;0>  & \stackrel{cube-edge}{\to} <100000000;1> \nonumber\\
  &  \stackrel{cube-edge}{\to} <110000000;1>.\nonumber
\end{align}
We next one hop along the ring to change the ring coordinate from $b=0$ to $b=4$:
$$<110000000;1>  \stackrel{ring-edge}{\to} <110000000;4>$$

Now with $b=4$, we flip $a_2$ and $a_3$ and then change
$b$ to $3$:
\begin{align}
 & <110000000;4> \stackrel{cube-edge}{\to} <110000100;4>\stackrel{cube-edge}{\to} \nonumber\\
& <110001100;4> \stackrel{ring-edge}{\to}<110001100;3>.\nonumber
\end{align}

Continue hopping alternatively across cube-edges and ring-edges, we get
\begin{align}
& <110001100;3> \stackrel{cube-edge}{\to}<110011100;3> \stackrel{cube-edge}{\to}\nonumber\\
& <110111100;3> \stackrel{ring-edge}{\to}
<110111100;2> \stackrel{cube-edge}{\to}\nonumber\\
& <111111100;2> \stackrel{ring-edge}{\to}<111111100;1> \stackrel{cube-edge}{\to}\nonumber \\
&<111111110;1>  \stackrel{cube-edge}{\to}<111111111;1>, \nonumber
\end{align}
at which point, we arrive at the same cube coordinate as node B, indicating that we are now in the same ring as node B. In this example, it then takes one more hope along the ring to get to node B:
$$<111111111;1>\stackrel{ring-edge}{\to}<111111111;2>.$$

It takes altogether $14$ hops to reach from the source node A to the
denotation node B, among which $9$ are cube-edge hops which change
the cube coordinate from $<000000000>$ to $<111111111>$, $4$ are
intermediate ring-edge hops which make the change of cube
coordinates possible, and $1$ is the final ring-edge hop to adjust
the ring coordinate. The total number of hops exceeds the
upper-bound provided in \cite{Hu2005} $k+j+1+\lfloor r/2 \rfloor =
2+7+1+2 = 12$.

\vspace{0.2cm}
{\bf \emph{Theorem 7: [Diameter of RCR]}}
\begin{itemize}
           \item If $\min_t (Num(k,r,j,t)) = 0,\ \ 0\leq t
\leq k+j-1$, then the RCR$(k,r,j)$ network is unconnected with a
diameter of $\infty$.
           \item If $\min_t (Num(k,r,j,t)) > 0,\ \ 0\leq t
\leq k+j-1$, then the diameter of RCR$(k,r,j)$ is upper bounded by
\begin{align}
    \label{equ:diameter of RCR}
     \left\{ \begin{array}{ll}
        k\!+\!j\!+\!r\!-\!1\!+\!\lceil\frac{r-1}{2}\rceil\!=\!k\!+\!j\!+\!r\!-\!1\!+\!\lfloor\frac{r}{2}\rfloor, & r \leq 3\\
        k\!+\!j\!+\!r\!+\!\lceil\frac{r-1}{2}\rceil\!-\!2\!=\!k\!+\!j\!+\!r\!-\!2\!+\!\lfloor\frac{r}{2}\rfloor, & r > 3
    \end{array}\right.
\end{align}
\end{itemize}
This bound is tight in both cases.

{\it Proof:} Consider an RCR$(k,r,j)$ network. A pair of nodes
 have the farthest distance when they
 stay in two rings whose cube coordinates differ in every bit
position. Without loss of generality, suppose the source node has
coordinate $<00\cdots 00;b=0>$ and the destination node has
coordinate $<11\cdots 11;b=b_0>$ for some valid value $b_0$.
From Example 6, the worst case involves many intermediate rings,
 such that one has to go
through all the possible values of $b$ in order to find
connecting cube edges
to reach the destination ring.
Depending on where $b_0$ is closer to $r-1$ or to $1$ in a ring, one
may choose to move clockwise or counter-clockwise along the
intermediate rings. With at the most $k+j$ cube-edge  hops and $r-1$
ring-edge hops, we will have arrived  either at $<11\cdots 11;r-1>$
or at $<11\cdots 11;1>$.

\begin{figure}[htbp]
        \centerline{
        \includegraphics[width=1.2in]{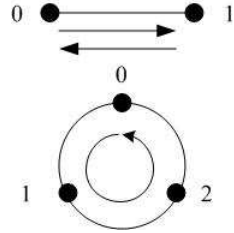}
}\vspace{-0.2cm}
        \caption{Long distance between a pair of nodes in RCR for $r \leq 3$.}
        \label{fig:diameter1}
    \end{figure}

\begin{figure}[htbp]
        \centerline{
        \includegraphics[width=1.2in]{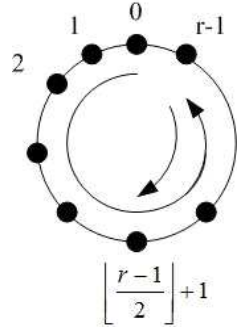}
}\vspace{-0.2cm}
        \caption{Long distance between a pair of nodes in RCR for $r>3$.}
        \label{fig:diameter}
    \end{figure}

Now to move along the destination ring to get to the destination
node, we have determine $b_0$ which has the longest distance with
$1$ and $r-1$. If $r \leq 3$, then $b_0=0$ as shown in
\ref{fig:diameter1}.  If $r >3$, then $1 < b_0 < r-1$. As shown in
figure \ref{fig:diameter}, the minimum distance between $b_0$ and
$r-1$ or $b_0$ and $1$ does not exceed
$\lfloor\frac{r}{2}\rfloor-1$. Therefore, it takes no more than a
total of $k+j+r+\lfloor\frac{r}{2}\rfloor-2$ hops to reach from any
node to any other node in an RCR$(k,r,j)$ network. It is easy to see
that this upper-bound is tight, since Example 6 achieves the bound
with equality.
\vspace{-0.5cm}
\vspace{-0.5cm}

\section{Modified RCR Networks}
\label{sec:modified RCR network}

We have thus far revisited RCR networks and rectified the results on
node degree, symmetry, connectivity, bisection width and network
diameter. A particularly desirable property of  RCR networks is
their easy construction and high scalability. However, the current
edge connecting rules have not fully exploited the potential
topological beauty of this class of networks. For example, Theorem 1
states that the node degree of an RCR$(k,r,j)$ network may not be
uniform, thus making network symmetry impossible. Further, Lemma 4
and Theorem 5 suggest that high-order RCR networks are doomed to be
unconnected, which significantly limits the ``useful'' scalability
of RCR networks. In what follows, we will modify the RCR networks in
\cite{Sun2000} by redefining the rule for cube edge connection. The
modified RCR networks, referred to as Class-II RCR networks, now
possess uniform node degree regardless of the parameters used, and
hence enjoy a better structural regularity and connectivity.

\vspace{0.2cm}
{\bf \emph{Definition 3:}} A Class-II recursive-cube-of-ring, denoted as
RCR-II$(k,r,j)$, is determined by three parameters, the cube
dimension $k\ge 0$, the ring dimension $r\ge 1$ and the level of
expansion $j\ge 0$. The construction of RCR-II$(k,r,j)$ is similar
to that of RCR$(k,r,j)$, where node coordinates are represented by:
\begin{align}
<\underbrace{a_{k+j-1},a_{k+j-2},...,a_{0}}_{cube\ coordinate};\underbrace{b}_{ring\ coordinate} > \nonumber\\
 \in \{0,1\}^{k+j}\times \{0,1,\cdots, r-1\},
\nonumber
\end{align}
and the nodes having the same cube coordinates but different ring
coordinates belong to the same ring. The only difference between
RCR-II and RCR is the connection of cube edges. In RCR-II, node
$<a_{k+j-1},\cdots, a_{t+1},a_{t}=0,a_{t-1},...,a_{0};b>$ can only
be connected to node  $<a_{k+j-1},\cdots,
a_{t+1},a_t=1,a_{t-1},...,a_{0};b>$ when the following constraints
are satisfied:
\begin{align}
t=g(bj+x,k+j)\stackrel{\Delta}{=}mod(bj+x,k+j),\\
0\le b \le   r-1, \ \ \ 
0 \leq x \leq k-1.
\end{align}

\vspace{0.2cm}
{\bf \emph{Theorem 8: [Node degree of RCR-II]}} An RCR-II$(k,r,j)$ network
has a uniform node degree $D_n=\min(k+r-1, k+2)$.

{\it Proof:} From linear algebraic, we know that $k$ distinct input
values for $x$ will yield $k$ distinct output values for $g(bj+x,
k+j)=mod(bj+x,k+j)$. Hence a node in RCR-II always has $k$ cube
neighbors. Since every node has $r-1$ ring neighbors for $r\le 2$
and $2$ ring neighbors for
 $r>2$, the result in Theorem 8 thus follows.

\begin{figure*}[htbp]
\centering
\includegraphics[width=4.5in]{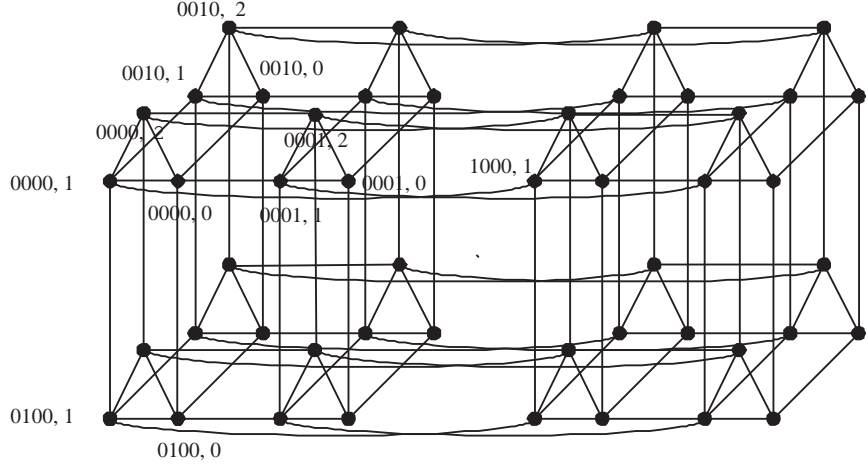}
\caption{Construction of RCR-II(3,3,1) network: uniform-node-degree
RCR-II network} \label{fig:nodedegree1}
\vspace{-0.3cm}
\end{figure*}

\vspace{0.2cm}
 {\bf \emph{Example 7: [Uniform node degree for RCR-II]}}
Example 1 shows that an RCR$(3,3,1)$ network has non-uniform node
degree. In comparison, RCR-II$(3,3,1)$, whose structure is depicted
in Fig.\ref{fig:nodedegree1}, has a uniform node degree of
$D_n=\min(3+3-1,3+2)=5$. To see this, note that via the definition
of RCR-II networks,
 $x$ may take  three possible values, $x\in\{0,1,2\}$, and the auxiliary variable $b$ may also take  three
possible values: $b\in\{0,1,2\}$. A node in this RCR-II network may
experience  one of the following three scenarios:
\begin{itemize}
  \item When $b\!=\!0$, the possible values for $f(bj \!+\!x, k\!+\!j)$
  are $0,1,2$, which correspond to $x\!=\!0,1,2$, respectively.
  The nodes in this case will have 3 cube neighbors,
  which lead to a node degree of $D_n=3+2=5$.
  \item When $b\!=\!1$, the possible values for $f(b j \!+\!x, k\!+\!j)$
  are $1,2,3$, so the nodes here also have 3 cube neighbors
  and a node degree of 5.
  \item When $b\!=\!2$, the possible values for $f(b j\! +\!x, k\!+\!j)$:
  are $2,3,0$, so the nodes here also have 3 cube neighbors
  and a node degree of 5
\end{itemize}
Therefore, the node degree is uniform as shown in Fig.
\ref{fig:nodedegree1}. However, this network is not symmetry. Below
 we discuss conditions for RCR-II to be
symmetric.

\vspace{0.2cm}
{\bf \emph{Theorem 9: [Symmetry of RCR-II]}}  An RCR-II$(k,r,j)$ network is symmetric if $mod(rj, k+j)=0$.

{\it Proof:} To show that RCR-II$(k,r,j)$ is symmetric, we need to
show that the network viewed from an arbitrary node
$<\alpha_{k+j-1},...,\alpha_0;\beta>$ has the same network topology
or neighbor-hood connectivity, as viewed from node $<00\cdots
00;0>$. This is equivalent to finding a proper rule that transforms
RCR-II$(k,r,j)$ to itself, such that the $<00\cdots 00;0>$  is
mapped to   $<\alpha_{k+j-1},...,\alpha_0;\beta>$,  and all the
other nodes and edges are mapped in a way that preserves the
original network topology and  cube- and ring-connecting rules.

Let the origin $<00\cdots 00;0>$ be  mapped to the new origin   $<\alpha_{k+j-1},...,\alpha_0;\beta>$. Assume that an arbitrary node $<a_{k+j-1},...,a_i,...,a_0;b>$ is correspondingly mapped to $<a'_{k+j-1},...,a'_i,...,a'_0;b'>$.
We define the transform as follows:
\begin{align}
    \label{equ:definition of mapping1}
     {\rm ring\ coordinate:\ } &   b'  =  mod(b + \beta, r); \\
     \label{equ:definition of mapping2}
    {\rm cube\ correlate:\ } &   a'_{mod(t-\beta k,k+j)}  =  a_t \oplus \alpha_{mod(t-\beta k,k+j)};
\end{align}
where the network parameters $k,r,j$ and the new origin $<\alpha_{k+j-1},\cdots,\alpha_0;\beta>$ are pre-determined constants.

(i) First, this transform is an enclosure, i.e.,  a valid node coordinate is mapped to a valid node coordinate. From (\ref{equ:definition of mapping1}),  the new ring coordinate $b'$ takes value between 0 and $r-1$ and is therefore a valid ring coordinate. From (\ref{equ:definition of mapping2}), the new cube coordinate $a'_{j}$ has index $0\le j \le k+j-1$ and takes value $a'_{j}\in\{0,1\}$, and
is therefore a valid cube coordinate.

(ii) Second, this transform is a one-to-one mapping. From (\ref{equ:definition of mapping1}) and (\ref{equ:definition of mapping2}), it is easy to see that if two node have different ring coordinates  and/or different cube coordinates before the transform, they will  take on different ring coordinates and/or different cube coordinates after transform.
Further, since the new ring coordinate is only a function of the old ring coordinate (and predetermined constants) and that the new cube coordinate is only a function of the old cube coordinate (and predetermined constants). If two node have the same ring- or cube-coordinates before the transform, they will take on the same ring- or cube-coordinates after transform.

(iii) Third, the ring edges are preserved under the transform. To see this, consider two nodes in the same ring. These nodes therefore have  a common cube coordinate but different ring coordinates. From (ii), after transform, they will take on a common cube coordinate and different ring coordinates, and hence are still in the same ring.

(iv) Finally, the cube edges are preserved under the transform.
Suppose there is a cube edge between $<a_{k+j-1},...,a_{t+1}, a_t,
a_{t-1}...,a_0;b>$ and
$<a_{k+j-1},...,a_{t+1},\bar{a}_t,a_{t-1}...,a_0;b>$, where,
according to the definition of RCR-II networks, $t = mod(b j + x_0,
k+j)$ for some integer value of  $x_0\in[0,k-1]$. Assume
$<a_{k+j-1},...,a_{t+1}, a_t, a_{t-1}...,a_0;b>$ is mapped to
$<a'_{k+j-1},...,a'_{t'+1}, a'_{t'}, a'_{t'-1}...,a'_0;b'>$ where  $
t'=mod(t-\beta k,k+j)$. From (\ref{equ:definition of mapping2}),
$<a_{k+j-1},...,a_{t+1}, \hat{a}_t, a_{t-1}...,a_0;b>$ is definitely
mapped to  $<a'_{k+j-1},...,a'_{t'+1}, \hat{a'}_{t'},
a'_{t'-1}...,a'_0;b'>$. Hence, it is sufficient to show that there
exists a cube edge connecting $<a'_{k+j-1},...,a'_{t'+1}, a'_{t'},
a'_{t'-1}...,a'_0;b'>$ and $<a'_{k+j-1},...,a'_{t'+1},
\hat{a'}_{t'}, a'_{t'-1}...,a'_0;b'>$, that is, $t'$ satisfies $t'=
mod (b'j+x, k+j)$ for some $x\in[0,k-1]$ (see the definition of
RCR-II networks).  We have

\begin{align}
t' &=  mod(t-\beta k,k+j),\\
    &= mod( (b j + x_0) - \beta k, k+j),\\
      &= mod( (b j+\beta j) +x_0 -(\beta k+\beta j), k+j),\\
      &= mod ((b+\beta)j+x_0, k+j). \label{equ:RCRII_symmetry_temp1}
\end{align}
From (\ref{equ:definition of mapping1}),  $b'=b+\beta +Ar$ for some
integer $A$.  From the assumption, $mod(rj,k+j) = 0$. Hence,
\begin{align}
&  mod((b+\beta)j+x_0,k+j) \nonumber\\
=& mod( (b +\beta)j + x_0 +Arj, k+j),\\
      =& mod( (b +\beta+Ar) j +x_0, k+j),\\
      =& mod (b'j+x_0, k+j). \label{equ:RCRII_symmetry_temp2}
\end{align}
Gathering (\ref{equ:RCRII_symmetry_temp1}) and (\ref{equ:RCRII_symmetry_temp2}), we get $t'=mod(b'j+x_0,k+j)$ where $0\le x_0\le k-1$. Hence, the cube connectivity is preserved after the transform.

\begin{figure}[htbp]
        \centering
        \includegraphics[width=2.2in,height=2.0in]{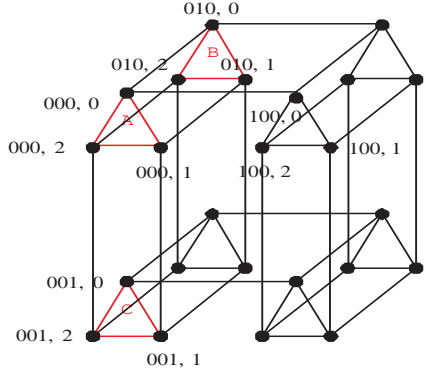}
        \begin{center}(A)asymmetric RCR(2,3,1)\end{center}
        \includegraphics[width=2.2in,height=2.0in]{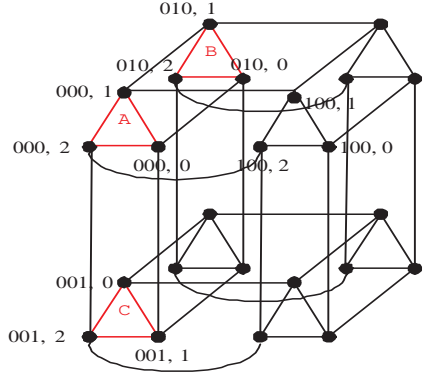}
        \begin{center}(B)symmetric RCR-II(2,3,1)\end{center}
        \caption{Construction of RCR(2,3,1) and RCR-II(2,3,1).}
        \label{fig:compare}
    \end{figure}

\vspace{0.2cm}
{\bf \emph{Example 8: [Symmetric RCR-II]}}
To help demonstrate the symmetry (and the balanced structure) of RCR-II networks, compare an RCR and RCR-II network with the same parameters
$(2,3,1)$ in Fig. \ref{fig:compare}. It is easy to see that
RCR(2,3,1) is asymmetric, since there does not exist a non-distorted
mapping that transforms node $<000,0>$ to node $<000,1>$. However, thanks to
the different edge connecting rules between ring $A$, $B$ and $C$, RCR-II(2,3,1) presents a symmetric network.

\vspace{0.2cm}
{\bf \emph{Theorem 10: [Sufficient and necessary condition for connectivity of RCR-II]}} An RCR-II$(k,r,j)$ network is connected if and only if
 $g(b j + x, k+j)=mod(bj+x,k+j)\covering [0, k+j-1]$ for $0\le b\le r-1$ and $0\le x\le k-1$.

{\it Proof:} The proof follows almost the same procedure as that for Theorem 4, and is therefore omitted.

\vspace{0.2cm}
{\bf \emph{Theorem 11: [Parameters for RCR-II to be connected]}}
An RCR-II$(k,r,j)$
network is connected, if and only if  $(r-1)k \geq j$.

{\it Proof}: We show that $(r-1)k\ge j$ is a necessary condition through proof-by-contradiction. Assume the network is connected but $(r-1) k <j$. Since $0\le b\le r-1$ and $ 0\le x \le k-1$,  there are $r$ possible values for $b$ and $k$ possible values for $x$, and hence at the most $rk$ different values for $t=\mod(bj+x, k+j)$. Since $rk<k+j$, it is impossible for $mod(bj+x,k+j)$ to cover $[0,k+j-1]$, so the network cannot be connected.

We now  show that $(r-1)k\ge j$ is also sufficient.

Case I: If $r=1$ and $(r-1) k \geq j$, then $j=0$ and $b=0$. It is easy to see that
 $g(bj + x, k+j)=mod(x,k)\covering [0, k-1]$ for $x\in[0,k-1]$. According to Theorem 10, the network is therefore connected.

Case II: If $r=2$ and $(r-1) k \geq j$, then $j \leq k$ and $b=0,1$.
We have $g(bj + x, k+j)=mod(x,k+j)\covering [0, k-1]$ when $b=0$ and $x\in[0,k-1]$,
 and $g(b j + x, k+j)=mod(j+x, k+j)\covering [j, j+k-1]$ when $b=1$ and $x\in[0,k-1]$. Since $j \leq k$, $g(bj+x,k+j)\covering [0, j+k-1]$ and the network is connected.

Case III: If $r>2$ and $(r-1)k \geq j$, we consider two subcases $k>j+1$ and $k\le j+1$.

(i) When $k \geq j$ and $x\in[0,k-1]$,
\begin{align}
& {\rm when\ } b=0,\\
& \ \ g(bj + x, k+j)=mod(x,k+j)\covering [0,k-1],\nonumber\\
& {\rm when\ } b=1, \\
& \ \ g(bj+x,k+j)=mod(j+x,k+j)\covering [j, k+j-1].\nonumber
\end{align}
Since $j \leq k$, $g(bj+x,k+j)\covering[0,k+j-1]$.

(ii) When $k \leq j-1$. Since $(r-1)k\geq j$, $j/k\le r-1$. Since
$j$, $k$ and $r$ are all integers, $\lceil  j/k \rceil \le r-1$.
Since $b\in[0,r-1]$, consider $b$ taking values from $0,1,\cdots,
\lceil j/k \rceil$.
\begin{align}
{\rm when \ } b=0,\ & g(bj+x, k+j) =mod(x,k+j)\nonumber\\
&\covering [0,k-1],\\
{\rm when \ }b=1,\ &g(bj+x,k+j) =mod(j+x,k+j)\nonumber\\
& \covering [j, j+k-1],\\
{\rm when \ }b=2, \ & g(bj+x,k+j)=mod(2j+x,k+j) \nonumber\\
& \covering [j-k,j-1], \\
{\rm when \ }b=3, \ & g(bj+x,k+j)=mod(3j+x,k+j)\nonumber\\
& \covering [j-2k ,j-k-1],\\
&  \cdots  \cdots\nonumber \\
 {\rm when \ }b\!=\!\lceil  j/k \rceil, & \  
g(bj\!+\!x,k\!+\!j)\!=\!mod(\lceil  j/k \rceil j\!+\!x,k\!+\!j) \nonumber\\
&  \ \covering
[j\!+\!k\!-\!\lceil j/k \rceil k,\ j\!+\!2k\!-\!\lceil j/k \rceil k\!-\!1],
\end{align}
Since $j+k- \lceil j/k \rceil  k  = k(j/k-\lceil j/k
\rceil+1)  \leq
 k$, 
we can see that $f(b j + x, k+j)\covering [0,k+j-1]$. Therefore, the
network is connected.

Class-II RCR networks exhibit similar topological properties for the
bisection width and the diameter as the original RCR networks.

\vspace{0.2cm}
{\bf \emph{Theorem 12: [Bisection Width of RCR]}} The bisection width of an RCR-II$(k,r,j)$ network is upper-bounded by:
\begin{equation}
B_{RCR-II}(k,r,j)\!\le\! \min_{t\in\{0,\cdots,k+j-1\}}\!\!\!\!  Num(k,r,j,t)\frac{N}{2r},
\end{equation}
where $r$ is the dimension of rings, $N=r2^{k+j}$ is the total number of nodes,  and $Num(k,r,j,t)$ the number of integer values $b\in [0,k-1]$ that satisfies $g(b j + x, k+j)=mod(bj+x,k+j) = t$ for given $k,r,j,t$, where $x\in[0,k-1]$.

\vspace{0.2cm}
{\bf \emph{Theorem 13: [Diameter of RCR-II]}}
\begin{itemize}
           \item If $\min_t (Num(k,r,j,t)) \!=\! 0,\  0\!\leq\! t
\!\leq\! k\!+\!j\!-\!1$, then the RCR-II$(k,r,j)$ network is unconnected with a
diameter of $\infty$.
           \item If $\min_t (Num(k,r,j,t))\! >\! 0,\ 0\!\leq\! t
\!\leq\! k\!+\!j\!-\!1$, then the diameter of RCR-II$(k,r,j)$ is upper bounded by
\begin{align}
    \label{equ:diameter of RCR-II}
     \left\{ \begin{array}{ll}
        k+j+r-1+\lceil\frac{r-1}{2}\rceil, & r \leq 3\\
        k+j+r+\lceil\frac{r-1}{2}\rceil-2, & r > 3
    \end{array}\right.
\end{align}
\end{itemize}
This bound is tight in both cases.

Theorems 12 and 13 can be proven using almost the identical arguments as those of RCR networks, and is therefore omitted.

\section{Conclusion}
Recursive-cube-of-ring (RCR) networks, proposed by Sun {\it et al}
\cite{Sun2000} and further analyzed by Hu {\it et al}
 \cite{Hu2005}, are a rich class of scalable interconnection networks that are determined
by three parameters, the ring dimension, the cube dimension and the
expansion number. Because of the many available combinations of
these three parameters, RCR networks take on a very rich pool of
possibilities with rather diverse structures, thus complicating the
analysis of their topological properties.

In this paper, we perform a close examination of RCR networks,
including the many special cases. Depending on the choice of the
parameters, RCR networks may expose rather different properties from
each other, some of which are less desirable for parallel computing.
For example, for the same seed network (i.e. the same ring diameter
and cube diameter), expanding an RCR network an additional level may
all of sudden change the network from well-connected to segmented.
Our contribution in the first part of this paper is the correction
of several misunderstanding and inaccuracies in the previous RCR
analysis,  including node degree, connectivity, symmetry, diameter
and bisection \cite{Sun2000,Hu2005}. Validating examples are
provided along with the discussion to support our analysis.

Since these RCR networks do not have uniform node degrees nor
possess network asymmetry, the second part of the paper focuses on
improving and enhancing this class of networks. Our contribution
here is the proposition of a class of modified RCR networks, termed
RCR-II networks, which preserve the same simplicity, richness and
scalability as the original RCRs, but which have uniform node
degrees irrespective of the network parameters, and exhibit better
connectivity and better symmetry than the original construction.
These better properties are achieved with only a simple change of
the cube edge connecting rules. Further, since uniform node degree
is but necessary condition for a network to be symmetric, sufficient
conditions to guarantee a symmetric RCR-II are also derived. Our
studies and findings in this paper provide a useful guidance for
choosing good parameters for RCR networks.


\end{document}